\documentclass[10pt]{amsart}
\usepackage{amsfonts}
\usepackage{amsmath}
\usepackage{amsthm}
\usepackage{amssymb}
\usepackage{latexsym}
\usepackage{multicol}
\usepackage{verbatim}
\usepackage{tabularx}
\usepackage{amscd}

\advance\textwidth by 1.2in \advance\oddsidemargin by -.6in
\advance\evensidemargin by -.6in
\newtheorem*{cor}{Corollary}
\newtheorem*{lem}{Lemma}
\newtheorem*{prop}{Proposition}

\theoremstyle{definition}
\newtheorem*{defn}{Definition}
\theoremstyle{definition}
\newtheorem*{thm}{Theorem}

\newtheorem*{rem}{Remark}

\newcounter{cnt}
\newenvironment{enumerit}{\begin{list}{{\hfill\rm(\roman{cnt})\hfill}}{%
\settowidth{\labelwidth}{{\rm(iv)}}\leftmargin=\labelwidth%
\advance\leftmargin by
\labelsep\rightmargin=0pt\usecounter{cnt}}}{\end{list}}

\theoremstyle{remark}

\numberwithin{equation}{section} 

\newtheorem*{example}{Example}
\newenvironment{ex}
 {\begin{example}\normalfont}
 {\end{example}}

\def\sop_#1^#2{\text{\scriptsize $\bigoplus\limits_{#1}^{#2}$}}

\newcommand{\nc}{\newcommand}
\newcommand{\rnc}{\renewcommand}

\nc{\cal}{\mathcal} \nc{\goth}{\mathfrak} \rnc{\bold}{\mathbf}

\makeatletter
\def\section{\def\@secnumfont{\mdseries}\@startsection{section}{1}%
  \z@{.7\linespacing\@plus\linespacing}{.5\linespacing}%
  {\normalfont\scshape\centering}}
\def\subsection{\def\@secnumfont{\bfseries}\@startsection{subsection}{2}%
  {\parindent}{.5\linespacing\@plus.7\linespacing}{-.5em}%
  {\normalfont\bfseries}}
\makeatother

\def\subl#1{\subsection{}\label{#1}}

\nc{\Cal}{\cal}  \nc{\ch}{\mbox{ch}} \nc{\Z}{{\bold Z}}
\nc{\J}{{\cal J}} \nc{\C}{{\bold C}} \nc{\Q}{{\bold Q}}

\nc{\N}{{\Bbb N}} \nc\boa{\bold a} \nc\bob{\bold b} \nc\boc{\bold
c} \nc\bod{\bold d} \nc\boe{\bold e} \nc\bof{\bold f}
\nc\bog{\bold g} \nc\boh{\bold h} \nc\boi{\bold i} \nc\boj{\bold
j} \nc\bok{\bold k} \nc\bol{\bold l} \nc\bom{\bold m}
\nc\bon{\bold n} \nc\boo{\bold o} \nc\bop{\bold p} \nc\boq{\bold
q} \nc\bor{\bold r} \nc\bos{\bold s} \nc\bou{\bold u}
\nc\bov{\bold v} \nc\bow{\bold w} \nc\boz{\bold z}

\nc\bb{\bold B} \nc\bc{\bold C} \nc\bd{\bold D}
\nc\be{\bold E}  \nc\bg{\bold G} \nc\bh{\bold H} \nc\bi{\bold I}
\nc\bj{\bold J} \nc\bk{\bold K} \nc\bl{\bold L} \nc\bm{\bold M}
\nc\bn{\bold N} \nc\bo{\bold O} \nc\bp{\bold P} \nc\bq{\bold Q}
\nc\br{\bold R} \nc\bs{\bold S} \nc\bt{\bold T} \nc\bu{\bold U}
\nc\bv{\bold V} \nc\bw{\bold W} \nc\bz{\bold Z} \nc\bx{\bold X}

\newcommand{\q}[2]{\genfrac{[}{]}{0pt}{}{#1}{#2}}

\newcommand{\ba}{\mathbb{A}}
\newcommand{\e}{\widetilde{e_i}}
\newcommand{\f}{\widetilde{f_i}}
\newcommand{\w}{\widetilde}
\newcommand{\qq}{\mathbb{Q}(q)}

\newcommand{\ns}{\subset\negthickspace\negthickspace\negthinspace\slash \ }

\begin{document}

\title{On Crystal Bases and Enright's Completions}
\author{Dijana Jakeli\'c}
\address{Department of Mathematics, University of Virginia, Charlottesville, VA 22904} \email{dj7m@virginia.edu}
\maketitle \setcounter{section}{0}

\flushbottom

\maketitle \centerline{\small{
\begin{minipage}{350pt}
{\bf Abstract:} We investigate the interplay of crystal bases and
completions in the sense of Enright on certain nonintegrable
representations of quantum groups. We define completions of crystal
bases, show that this notion of completion is compatible with
Enright's completion of modules, prove that every module in our
category has a crystal basis which can be completed and that a
completion of the crystal lattice is unique. Furthermore, we give
two constructions of the completion of a crystal lattice.
\end{minipage}}}

\bigskip

\hspace{3.4in} \emph{Dedicated to the memory of Des Sheiham}
\medskip
\section*{Introduction}

In this paper we are concerned with bringing together two distinct
notions of the representation theory of quantized enveloping
algebras -- those of crystal bases and completions.

The theory of crystal and canonical bases is one of the most
remarkable developments in representation theory. It was introduced
independently by M. Kashiwara \cite{K1,K2} and G. Lusztig \cite{L}
in a combinatorial and a geometric way, respectively. In this paper
we follow Kashiwara's approach. Roughly speaking, a crystal basis of
an integrable representation for the quantized enveloping algebra
$U_q(\mathfrak{g})$ of a finite-dimensional complex semisimple Lie
algebra $\mathfrak{g}$ (or more generally, of a symmetrizable
Kac-Moody Lie algebra) is a pair consisting of a lattice of the
module, called a crystal lattice, and a vector space basis of a
quotient of the crystal lattice, called a crystal. It is actually a
certain parametrization of bases of the module with a number of
desirable properties and it encodes the intrinsic structure of the
module in a combinatorial way. One of the most important
combinatorial realizations of crystals is Littelmann's path model
\cite {Li1, Li2}.

T. Enright \cite{E} introduced a notion of completion with respect
to a simple root of $\mathfrak{g}$ on the category
$\mathcal{I}(\mathfrak{g})$ of weight modules of $\mathfrak{g}$ that
are $U^-(\mathfrak{g})$-torsion free and $U^+(\mathfrak{g})$-finite.
Completion is an effective process of obtaining new representations
from a given one, containing the original one as a
subrepresentation. Enright used the completion functor in order to
algebraically construct the fundamental series representations. Soon
after, V. Deodhar \cite{D} realized the completion functor via Ore
localization giving a concrete model of completion and used it to
prove Enright's uniqueness conjecture arising from considering
successive completions (stressing the ``$\mathfrak{sl_2}$-nature''
of completions). A. Joseph \cite{Jo} then generalized this functor
to the Bernstein-Gelfand-Gelfand category
$\mathcal{O}(\mathfrak{g})$ and gave a refinement of the Jantzen
conjecture by studying the lifting of a contravariant form under the
action of the completion functor. Later, Y.M. Zou \cite{Z} extended
some of the results obtained by Enright and Deodhar to the quantum
groups setting.

The starting motivation of this work was to study both crystal bases
and completions of modules belonging to the category
$\mathcal{I}\left(U_q(\mathfrak{g})\right)$ (thus nonintegrable
ones) and examine how the two concepts relate to each other with the
aim of introducing a notion of completion of crystal bases which
would be compatible with Enright's completion of modules. It is
natural to expect that such a program may eventually lead to some
potentially interesting interplay between the theory of crystal
bases and that of fundamental series representations. Furthermore,
since all the Verma modules in a lattice of inclusions of Verma
modules may be obtained from the corresponding irreducible one by
means of completions, one may expect that a successful crystal base
theory related to completions could produce a combinatorial tool
relevant for studying Jantzen filtrations and some other
Kazhdan-Lusztig Theory related topics.

Beside defining crystal bases of integrable representations,
Kashiwara \cite{K2} also defined the crystal basis of the
quantization $U_q^-(\mathfrak{g})$ of the universal enveloping
algebra of the nilpotent part of $\mathfrak{g}$ by considering
$U_q^-(\mathfrak{g})$ as a module for Kashiwara's algebra
$B_q(\mathfrak{g})$. Although the former was the main goal of
\cite{K2}, the latter provided a way to simultaneously consider the
bases of all integrable representations, and the interplay of the
two notions of crystal bases played a central role in the paper.
Since \cite{K2} there have been only a few papers concerned with
crystal bases or crystals of representations which are not
necessarily integrable, including \cite{KKM}, \cite{NZ}, and
\cite{Ze}. Also, in a joint work with V. Chari and A. Moura
\cite{CJM} we considered the problem of tensor product
decompositions into indecomposables for nonintegrable modules in the
BGG category $\mathcal{O}$ by introducing the combinatorial objects
called branched crystals which satisfy a relaxed axiom for formal
invertibility of Kashiwara's operators.

In this paper, we follow Kashiwara's definition of crystal bases,
thus the $U_q^-$-torsion free modules in question are naturally
viewed as $\mathcal{B}_q(\mathfrak{g})$-modules. On the other hand,
in regard to completions they are thought of as
$U_q(\mathfrak{g})$-modules. This aspect makes the situation more
interesting and a synchronization of the two structures becomes
essential. As the $U_q(\mathfrak{sl_2})$-case is already quite
intricate, we restrict ourselves to that case in this paper leaving
the consideration on how to extend this theory to the higher rank
case to a future work.

We introduce a category $\w{\mathcal{I}}$ consisting of
$U_q$-modules in
$\mathcal{I}=\mathcal{I}\left(U_q(\mathfrak{sl_2})\right)$ with a
compatible $\mathcal{B}_q$-structure and obtain a decomposition of
every module in $\w{\mathcal{I}}$ into a direct sum of
indecomposable $U_q$-submodules which are also
$\mathcal{B}_q$-invariant (cf. Theorem \ref{t:UB}). Due to this
decomposition, we are able to introduce weight spaces into our
consideration of crystal bases of modules in $\w{\mathcal{I}}$ even
though these crystal bases arise from the
$\mathcal{B}_q$-structures. Therefore, we get a setting resembling
the one of integrable $U_q$-modules. We define a notion of a
complete crystal basis and a completion of a crystal basis. In the
process, we take advantage of a symmetric action of the Kashiwara
operators $\w{e}$ and $\w{f}$ on crystal bases. However, unlike the
case of modules where the completion of a module contains the
module, one cannot expect a completion of a crystal lattice to
contain the crystal lattice; this being indeed clear from the case
of Verma modules. Our definition of completion of bases involves a
very natural connection of crystal bases arising from
$\mathcal{B}_q$-structures with the ones arising from
$U_q$-structures. We prove that a crystal basis is complete if and
only if its corresponding module is. We also show that every module
in $\w{\mathcal{I}}$ has a crystal basis which can be completed and
moreover a completion of the crystal lattice is unique (cf. Theorems
\ref{t:CLC} and \ref{v:down}.).

The paper is arranged as follows. In Section 1, we set the notation
and review relevant results on completions and crystal bases. We
think it instructive to have this section done for
$U_q(\mathfrak{g})$, where $\mathfrak{g}$ is any simple Lie algebra,
as it is then more transparent how crucial a step the
$U_q(\mathfrak{sl_2})$-case is in solving the problem. However, the
remaining sections treat the $U_q(\mathfrak{sl_2})$-case, except
where stated otherwise. In Section 2, we consider some natural
$\mathcal{B}_q$-structures on the indecomposable $U_q$-modules in
$\mathcal{I}$, collect the desirable properties that a
$\mathcal{B}_q$-structure should have with respect to a
$U_q$-structure in order to define the category $\w{\mathcal{I}}$,
and prove the simultaneous decomposition theorem for modules in
$\w{\mathcal{I}}$. In Section 3, we look into the crystal bases with
which modules in $\w{\mathcal{I}}$ are naturally endowed via their
$\mathcal{B}_q$-structures, define complete crystal bases, show they
correspond to complete modules, define a completion of a crystal
basis, and prove the aforementioned Theorems \ref{t:CLC} and
\ref{v:down}. In Section 4 we give two constructions of the
completion of a crystal lattice -- one obtained by modifying
Deodhar's model of completion of modules and the other by applying
an operator used by Kashiwara in \cite{K1} to construct the
operators $\w{e}$ and $\w{f}$.

\bigskip

%%%%%%%%%%%%%%%%%%%%%%%%%%%%%%%%%%%%%%%%%%%%%%%%%%%%%%%%

\section{Preliminaries}

\subl{} Let $(a_{ij})_{1 \le i,j \le l}$ be the Cartan matrix of a
finite-dimensional complex simple Lie algebra $\mathfrak{g}$ and
$d_i$ unique positive integers such that $\mbox{gcd} \ (d_1, \dots,
d_l)=1$ and the matrix $(d_ia_{ij})_{1 \le i,j \le l}$ is symmetric.
Let $q$ be an indeterminate and $\qq$ the field of rational
functions of $q$ with coefficients in $\mathbb{Q}$. Set
$q_i=q^{d_i}$. The quantized enveloping algebra $U_q(\mathfrak{g})$
is the $\qq$-algebra with generators $e_i, f_i, t_i, t_i^{-1}, \ 1
\le i \le l$, and defining relations
\begin{align*}
 &t_it_i^{-1}=1=t_i^{-1}t_i, \ t_it_j=t_jt_i,\\
 &t_ie_jt_i^{-1}=q_i^{a_{ij}}e_j,  \ t_if_jt_i^{-1}=q_i^{-a_{ij}}f_j,\\
 &e_if_j-f_je_i=\delta_{ij} \frac{t_i-t_i^{-1}}{q_i-q_i^{-1}},\\
 &\sum_{s=0}^{1-a_{ij}} (-1)^s e_i^{(s)}e_je_i^{(1-a_{ij}-s)}=0, \ \sum_{s=0}^{1-a_{ij}} (-1)^s f_i^{(s)}
  f_jf_i^{(1-a_{ij}-s)}=0 \ \ (i \ne j)
\end{align*}
where $\displaystyle{e_i^{(n)}=\frac{e_i^n}{[n]_{i}!}, \ f_i^{(n)}=\frac{f_i^n}{[n]_{i}!}}$,
$[n]_i ! = [1]_i [2]_i \dots [n]_i \ (n \in \mathbb{Z^+})$, and
$\displaystyle{[n]_i = \frac{q_i^n-q_i^{-n}}{q_i-q_i^{-1}} \ (n \in \mathbb{Z})}$.

\medskip

Denote by $U_q^+(\mathfrak{g})$ (respectively $U_q^-(\mathfrak{g})$)
the subalgebra of $U_q(\mathfrak{g})$ generated by $e_i$
(respectively $f_i$), $1 \le i \le l$. Let $U_q^0(\mathfrak{g})$ be
the subalgebra of $U_q(\mathfrak{g})$ generated by $t_i$ and
$t_i^{-1}, 1 \le i \le l$. Multiplication defines an isomorphism of
$\qq$-vector spaces $U_q^-(\mathfrak{g})\otimes
U_q^0(\mathfrak{g})\otimes U_q^+(\mathfrak{g}) \to
U_q(\mathfrak{g})$.

Denote by ${U_q(\mathfrak{g})}_i$ the subalgebra of
$U_q(\mathfrak{g})$ generated by $e_i, f_i,t_i$ and $t_i^{-1}$.
Then ${U_q(\mathfrak{g})}_i \cong U_{q_i}(\mathfrak{sl_2})$.

\subl{}
Let $\mathfrak{h}$ be a Cartan subalgebra of $\mathfrak{g}$ and $\Phi$ the root system of
$(\mathfrak{g}, \mathfrak{h})$.
Let $\{\alpha_i\}_{1 \le i \le l}\subset \mathfrak{h}^*$ be a set of simple roots of $\mathfrak{g}$
and $\{h_i\}_{1 \le i \le l} \subset \mathfrak{h}$ the corresponding set of coroots so that $\alpha_j(h_i)=a_{ij}$. Let $Q=\sum_{i=1}^l \mathbb{Z} \alpha_i $ be the
root lattice and $P=\{\lambda \in \mathfrak{h}^* | \ \lambda(h_i) \in \mathbb{Z}, 1 \le i \le l\}$
the weight lattice
for $\mathfrak{g}$. Also, set $Q^+=\sum_{i=1}^l \mathbb{Z}^+ \alpha_i $ and
$P^+=\{\lambda \in \mathfrak{h}^* | \ \lambda(h_i) \in \mathbb{Z}^+, 1 \le i \le l\}$.
Denote by $W$ be the Weyl group of $\Phi$ and by $s_{i}$ the reflection with respect
to the simple root $\alpha_{i}$ $(1 \le i \le l)$.

For a $U_q(\mathfrak{g})$-module $M$ and $\lambda \in P$, we define the $\lambda$-weight space of $M$ as $M_\lambda=\{m \in M | \ t_im=q_i^{\lambda(h_i)}m, 1 \le i \le l\}$. $M$ is a weight module if it is a direct sum of its weight spaces.
If in addition there exist $\lambda \in P$ and a nonzero vector $m \in M_\lambda$
such that $e_i \cdot m=0$ for all $i$ and
$M=U_q(\mathfrak{g}) \cdot m$, then $M$ is a highest weight module with highest weight $\lambda$ and highest weight vector $m$.

The Verma module $M(\lambda)$ is the $U_q(\mathfrak{g})$-module
$U_q(\mathfrak{g}) \slash J(\lambda)$ where $J(\lambda)$ is the left
ideal of $U_q(\mathfrak{g})$ generated by $\{t_i-q_i^{\lambda(h_i)},
e_i | \ 1 \le i\le l\}$. $M(\lambda)$ is the universal highest
weight module of weight $\lambda$. The unique irreducible quotient
$V(\lambda)$ of $M(\lambda)$ is finite-dimensional iff $\lambda \in
P^+$.

%%%%%%%%%%%%%%%%%%%%%%%%%%%%%%%%%%%%%%%%%%%%%%%%%%%%%%%

\subl{}\label{cati} We recall the definitions and some properties of category $\mathcal{I}$, completions,
and $T$ modules
(cf. \cite{D}, \cite{E}, \cite{J}, \cite{Z}).

Let $\mathcal{I}\left(U_q(\mathfrak{g})\right)$ be the category of
$U_q(\mathfrak{g})$-modules $M$
satisfying:
(i) $M$ is a weight module,
(ii) $U_q^-(\mathfrak{g})$-action on $M$ is torsion free, and
(iii) $M$ is $U_q^+(\mathfrak{g})$-finite, i.e., $e_i$ acts locally nilpotently
on $M$ for all $i$.

Fix $i \in \{1, \dots, l\}$. Set
$M_n = \{m \in M| \ t_im=q_i^n m\} \ \mbox{for} \ n \in \mathbb{Z}, \
M^{e_i} = \{m \in M| \ e_im=0\} , \ \mbox{and} \
M_n^{e_i} = M_n \cap M^{e_i}$.
A module $M$ in $\mathcal{I}\left( U_q(\mathfrak{g}) \right)$ is
said to be complete with respect
to $i$ if $ f_i^{n+1}: \ M_n^{e_i} \to  M_{-n-2}^{e_i}$
is bijective for all $n \in \mathbb{Z}^+$.
A module $N$ in $\mathcal{I}\left( U_q(\mathfrak{g}) \right)$ is
a completion of
$M$ with respect to $i$ provided:
(i) $N$ is complete with respect to $i$,
(ii) $M$ is imbedded in $N$, and
(iii) $N \slash M$ is $f_i$-finite.

\begin{thm}(cf. \cite{D}, \cite{E}, \cite{Z})
\begin{enumerit}
\item Every module $M$ in $\mathcal{I}\left( U_q(\mathfrak{g}) \right)$ has a completion $C_i(M)$ with respect to $i$ $(1 \le i \le l)$ and any two such completions are naturally isomorphic.
\item Let $w \in W$ and $M \in \mathcal{I}\left(U_q(\mathfrak{g})\right)$. For any two reduced expressions
$w=s_{i_1} \dots s_{i_k}=s_{j_1} \dots s_{j_k}$,
there exists an isomorphism
$F: C_{i_1}\left(C_{i_2} \dots \left(C_{i_k}(M)\right) \dots \right) \to C_{j_1}
 \left(C_{j_2} \dots \left(C_{j_k}(M)\right) \dots \right)$ such that $F|_M$ is the identity.
\end{enumerit}
\end{thm}

Thus the process of completion depends essentially on the
$U_q(\mathfrak{sl_2})$-representation theory and so we pay special
attention to the case $\mathfrak{g}=\mathfrak{sl_2}$.

\subl{}
For brevity, write $U_q=U_q(\mathfrak{sl_2})$ and
$\mathcal{I}=\mathcal{I}\left(U_q(\mathfrak{sl_2})\right)$. Denote the generators of $\mathfrak{sl_2}$ by $e,f,t^{\pm 1}$.

Let $n \in \mathbb{Z}$. The quantum Casimir element $\displaystyle{C
= \frac{qt+q^{-1}t^{-1}}{(q-q^{-1})^2}} + fe$ acts on the Verma
module $M(n)$ as multiplication by scalar
$\displaystyle{c_n=\frac{q^{n+1}+q^{-n-1}}{(q-q^{-1})^2}}$. Consider
the left ideal $I(n)$ of $U_q$ defined by $ I(n)=U_q \{t- q^{-n-2},
\, e^{n+2}, \, (C-c_n)^2\}$. The $T$ module $T(n)$ is defined as
$T(n)=U_q \slash {I(n)}$. It is an indecomposable module belonging
to the category $\mathcal{I}$ and $0 \to M(n) \to T(n) \to M(-n-2)
\to 0$ is exact.

\begin{thm} \label{t:e}(cf. \cite{E}, \cite{Z})
\begin{enumerit}
\item The $M(n) \ (n \in \mathbb{Z})$ and the $T(n) \ (n \in \mathbb{Z}^+)$
are precisely all the indecomposable objects of the
category $\mathcal{I}$. Among these, the $M(n)$ for $n\ge -1$ and
$T(n)$ for $n \in \mathbb{Z}^+$ are the complete ones.
The completion of $M(-n-2)$ is $M(n)$ for $n \ge -1$.
\item Every module in $\mathcal{I}$ is a direct sum (not necessarily finite) of indecomposable
ones.
\end{enumerit}
\end{thm}

%%%%%%%%%%%%%%%%%%%%%%%%%%%%%%%%%%%%%%%%%%%%%%%%%%%%%

\subl{}
In the subsections that follow, we recall the main results on Kashiwara's crystal bases. (cf. \cite{K2}.)

Let $M$ be a finite-dimensional $U_q(\mathfrak{g})$-module. Fix an
index $i \ (1 \le i \le l)$. By the representation theory of
$U_q(\mathfrak{sl_2})$, $ \ M=\bigoplus_{\substack{\lambda \in P\\ 0
\le k \le \lambda(h_i)}} f_i^{(k)} \left(\textnormal{Ker} \, e_i
\cap M_\lambda\right)$. Hence, every element $u \in M_\lambda$ can
be uniquely written as $\displaystyle{u=\sum_{k \ge 0}
f_i^{(k)}u_k}$ where $u_k \in \textnormal{Ker} \, e_i \cap
M_{\lambda+k\alpha_i}$. The Kashiwara operators $\widetilde{e_i}$
and $\widetilde{f_i}$ are defined by
$\displaystyle{\widetilde{e_i}u=\sum_{k \ge 1} f_i^{(k-1)} u_k \quad
\textnormal{and} \quad \widetilde{f_i}u=\sum_{k \ge 0} f_i^{(k+1)}
u_k}$. It follows that $\widetilde{e_i}M_\lambda \subseteq
M_{\lambda+\alpha_i} \quad \textnormal{and} \quad
\widetilde{f_i}M_\lambda \subseteq M_{\lambda-\alpha_i}$. Also, if
$u \in \textnormal{Ker} \, e_i$ and $n>0$, then
$\widetilde{f_i}^nu=f_i^{(n)}u$.

Let $\mathbb{A}$ be the subring of $\mathbb{Q}(q)$ consisting of rational functions regular at $q=0$.
\begin{defn}\label{d:CLCB}
(\cite{K2}) \
A free $\ba$-submodule $L$ of $M$ is called a crystal lattice if:
(a) \ $M \cong \mathbb{Q}(q) \otimes_{\ba} L$,
(b) \ $L=\oplus_{\lambda \in P} L_\lambda$ \ where \ $L_\lambda=L \cap M_\lambda$, and
(c) \ $\widetilde{e_i}L \subseteq L$ \ and \ $\widetilde{f_i}L \subseteq L, \ 1 \le i \le l$.
A crystal basis of $M$ is a pair $(L,B)$ satisfying the following conditions:
(i) \ $L$ is a crystal lattice of $M$,
(ii) \ $B$ is a $\mathbb{Q}$-basis of $L \slash qL$,
(iii) \ $B=\sqcup_{\lambda \in P} B_\lambda$ where $B_\lambda=B \cap \left(L_\lambda
\slash qL_\lambda\right)$,
(iv) \ $\widetilde{e_i}B \subseteq B \cup \{0\}$ and $\widetilde{f_i}B \subseteq B \cup \{0\},
1 \le i \le l$, and
(v) \ for $b, b' \in B, \widetilde{f_i}b=b'$ if and only if $b=\widetilde{e_i}b'$.
\end{defn}

For $\lambda \in P^+$, consider the finite-dimensional irreducible
$U_q(\mathfrak{g})$-module $V(\lambda)$ with highest weight
$\lambda$ and highest weight vector $u_\lambda$. Let $L(\lambda)$ be
the smallest $\ba$-submodule of $V(\lambda)$ containing $u_\lambda$
which is stable under $\widetilde{f_i}$'s, i.e., $L(\lambda)$ is the
$\ba$-span of the vectors of the form $\widetilde{f_{i_1}} \dots
\widetilde{f_{i_r}} u_\lambda$, where $1 \le i_j \le l$ and $r \in
\mathbb{Z}^+$. Set $B(\lambda)=\{b \in L(\lambda) \slash qL(\lambda)
| \ b=\widetilde{f_{i_1}} \dots \widetilde{f_{i_r}} u_\lambda \
\textnormal{mod} \ qL(\lambda)\} \smallsetminus \{0\}$. Then
$\left(L(\lambda), B(\lambda)\right)$ is a crystal basis of
$V(\lambda)$ (cf. \cite{K2}, Theorem 2).

Since crystal bases are stable under direct sums, every
finite-dimensional $U_q(\mathfrak{g})$-module $M$ has a crystal
basis.

\subl{}
Isomorphism of crystal bases is defined as follows:

\begin{defn}\label{d:iso}
Let $(L_1,B_1)$ and $(L_2,B_2)$ be crystal bases of finite
dimensional $U_q(\mathfrak{g})$-modules $M_1$ and $M_2$,
respectively. We say that $(L_1,B_1) \cong (L_2,B_2)$ if there
exists a $U_q(\mathfrak{g})$-isomorphism $\varphi: M_1 \to M_2$
which induces an isomorphism of $\ba$-lattices $\varphi: L_1 \to
L_2$ such that $\bar{\varphi}(B_1) = B_2$ where $\bar{\varphi}:
L_1/qL_1 \to L_2/qL_2$ is the induced isomorphism of
$\mathbb{Q}$-vector spaces.

\end{defn}

 If $(L,B)$ is a crystal basis of a finite-dimensional $U_q(\mathfrak{g})$-module $M$,
then there exists an isomorphism $M \cong \oplus_{j} V(\lambda_j)$
by which $(L,B)$ is isomorphic to $\oplus_{j} \left(L(\lambda_j),
B(\lambda_j)\right)$ (cf. \cite{K2}, Theorem 3).

\subl{} Kashiwara's algebra $\mathcal{B}_q(\mathfrak{g})$ (cf.
\cite{K2}, \cite{N}) is the $\mathbb{Q}(q)$-algebra generated by
$e_i'$ and $f_i, \ 1 \le i \le l$, subject to relations
\begin{align}
 &e_i'f_j=q_i^{-a_{ij}}f_je_i'+\delta_{ij}\\
 &\sum_{s=0}^{1-a_{ij}} (-1)^s (e_i')^{(s)}(e_j')(e_i')^{(1-a_{ij}-s)}=0, \ \sum_{s=0}^{1-a_{ij}}
 (-1)^s f_i^{(s)}f_jf_i^{(1-a_{ij}-s)}=0 \ \ (i \ne j).
\end{align}

For the case $\mathfrak{g}=\mathfrak{sl_2}$, write $\mathcal{B}_q=\mathcal{B}_q(\mathfrak{sl_2})$.
Then the generators $e'$ and $f$ of $\mathcal{B}_q$ satisfy
\begin{equation}\label{b4}
 e'f=q^{-2}fe'+1.
\end{equation}

For each $P \in U_q^-(\mathfrak{g})$, there exist unique $R, Q \in
U_q^-(\mathfrak{g})$ such that $[e_i,
P]=\displaystyle{\frac{t_iQ-t_i^{-1}R}{q_i-q_i^{-1}}}$. Define $e_i'
\in \textnormal{End} \left(U_q^-(\mathfrak{g})\right)$ by
$e_i'(P)=R$. Then, $U_q^-(\mathfrak{g})$ is a left
$\mathcal{B}_q(\mathfrak{g})$-module where $e_i'$ acts as described
and $f_i$ acts by the left multiplication. Moreover,
$U_q^-(\mathfrak{g})$ is an irreducible
$\mathcal{B}_q(\mathfrak{g})$-module and $U_q^-(\mathfrak{g}) \in
\mathcal{O}\left(\mathcal{B}_q(\mathfrak{g})\right)$, where
$\mathcal{O}\left(\mathcal{B}_q(\mathfrak{g})\right)$ is the
category of $\mathcal{B}_q(\mathfrak{g})$-modules $M$ satisfying
that for all $u \in M$, there exists $r \ge 0$ such that
$e_{i_1}'e_{i_2}' \dots e_{i_r}' u=0$ for any $1 \le i_1, \dots, i_r
\le l$. Furthermore,
$\mathcal{O}\left(\mathcal{B}_q(\mathfrak{g})\right)$ is semisimple
and $U_q^-(\mathfrak{g})$ is the unique irreducible object of
$\mathcal{O}\left(\mathcal{B}_q(\mathfrak{g})\right)$, up to
isomorphism.
\begin{rem}\label{r:e}
For each $i$,
$e_i'(1)=0$ and
$e_i'(f_i^p)=q_i^{-(p-1)}[p]_if_i^{p-1}$ for $p \ge 1$.
\end{rem}

\subl{}
Let $M$ belong to the category
$\mathcal{O}\left(\mathcal{B}_q(\mathfrak{g})\right)$.
Fix an index $i (1 \le i \le l)$. Then
\begin{equation}\label{cb6}
 M=\oplus_{k \ge 0} f_i^{(k)} \ \mbox{Ker} \, e_i' \ .
\end{equation}
Define Kashiwara's operators $\widetilde{e_i}, \widetilde{f_i} \in \ \mbox{End} \, (M)$ by
\begin{equation}\label{cb7}
 \widetilde{e_i}\left(f_i^{(k)}u\right)=
                                        \begin{cases}
                                          f_i^{(k-1)}u,    &k \ge 1\\
                                          0,               &k=0
                                        \end{cases}
 \quad \quad \quad \quad
 \widetilde{f_i}\left(f_i^{(k)}u\right)=f_i^{(k+1)}u,
\end{equation}
for $u \in \ \mbox{Ker} \, e_i'$ and extend linearly.

Then $\widetilde{e_i} \widetilde{f_i}=1$. Also, $\widetilde{f_i} \widetilde{e_i}$ is the
projection onto $f_i M$ with respect
to $M=f_iM \oplus \ \mbox{Ker} \, e_i' \ $.
\begin{defn}\label{d:CLBCBB}
(\cite{K2}) \
A free $\ba$-submodule $L$ of $M$ is called a crystal lattice if:
(a) \ $M \cong \mathbb{Q}(q) \otimes_\ba L$, and
(b) \ $\widetilde{e_i}L \subseteq L$ and $\widetilde{f_i}L \subseteq L, 1 \le i \le l$.
A crystal basis of $M$ is a pair $(L,B)$ satisfying the following conditions:
(i) \ $L$ is a crystal lattice of $M$,
(ii) \ $B$ is a $\mathbb{Q}$-basis of $L \slash qL$,
(iii) \ $\widetilde{e_i}B \subseteq B \cup \{0\}$ and $\widetilde{f_i}B \subseteq B, 1 \le i \le l$, and
(iv) \ if $b \in B$ such that $\widetilde{e_i}b \in B$, then $b=\widetilde{f_i}\widetilde{e_i}b$.
\end{defn}
Isomorphism of crystal bases of modules in $\mathcal{O}\left(\mathcal{B}_q(\mathfrak{g})\right)$ is defined
analogously to isomorphism of crystal bases of finite-dimensional modules replacing $U_q(\mathfrak g)$ by $\mathcal B_q(\mathfrak g)$ in Definition \ref{d:iso}.

It will be clear from the context if we mean crystal basis in the sense of Definition \ref{d:CLCB} or in the
sense of Definition \ref{d:CLBCBB}. Otherwise, we will make the distinction.

Let $L(\infty)$ be the smallest $\ba$-submodule of
$U_q^-(\mathfrak{g})$ containing 1 that is stable under $\f$'s,
i.e., $L(\infty)=\ba$-span of $\{\widetilde{f_{i_1}} \dots
\widetilde{f_{i_r}} \cdot 1 | \ 1 \le i_j \le l; r \in
\mathbb{Z}^+\}$. Set $B(\infty)=\{b \in L(\infty) \slash qL(\infty)
| \ b=\widetilde{f_{i_1}} \dots \widetilde{f_{i_r}} \cdot 1 \,
\mbox{mod} \, qL(\infty)\}$. Then $\left(L(\infty),
B(\infty)\right)$ is a crystal basis of $U_q^-(\mathfrak{g})$
(\cite{K2}, Theorem 4). Moreover, any crystal basis of
$U_q^-(\mathfrak{g})$ coincides with $\left(L(\infty),
B(\infty)\right)$ up to a constant multiple.

The relation of $\left(L(\infty), B(\infty)\right)$ to
$\left(L(\lambda), B(\lambda)\right)$ is given by \cite{K2}, Theorem
5.

%%%%%%%%%%%%%%%%%%%%%%%%%%%%%%%%%%%%%%%%%%%%%%%%%%%%%%%%%%%%%%%%%%%%%%%

\section{$U_q$-modules with $\mathcal{B}_q$-module Structure}

In this section we are concerned with $U_q$-modules endowed with
$B_q$-module structures which allow decompositions that are
simultaneously $U_q$ and $B_q$-invariant. In \ref{bverma} we drop
the assumption $\mathfrak{g}=\mathfrak{sl_2}$, but we uphold it
otherwise.

\subl{}\label{bverma}
For $\lambda \in P$,
we consider $M(\lambda)$, the Verma module
with heighest weight $\lambda$.
If $m_\lambda$ is a highest weight vector of $M(\lambda)$,
then $M(\lambda)=U_q(\mathfrak{g})m_\lambda=U_q^-(\mathfrak{g})m_\lambda$.
Since $U_q^-(\mathfrak{g})$ is a $\mathcal{B}_q(\mathfrak{g})$-module, $M(\lambda)$ has a natural
$\mathcal{B}_q(\mathfrak{g})$-structure
via composition of $\mathbb{Q}(q)$-algebra homomorphisms
\[
 \begin{CD}
  \mathcal{B}_q(\mathfrak{g}) @>>> \textnormal{End} \left(U_q^-(\mathfrak{g})\right) @>\Theta>>
  \textnormal{End} \left(M(\lambda)\right).
 \end{CD}
\]
\noindent Namely, if $\varphi_\lambda : U_q^-(\mathfrak{g}) \to M(\lambda)$ is the
$U_q^-(\mathfrak{g})$-module
isomorphism sending 1 to $m_\lambda$, then $\Theta(g)=\varphi_\lambda g\varphi_\lambda^{-1}$
for $g \in \textnormal{End} \left(U_q^-(\mathfrak{g})\right)$.

Hence, if $f_i, e_i' \in \textnormal{End} \left(U_q^-(\mathfrak{g})\right), 1 \le i \le l$,
are defined as in \ref{r:e}, then
$f_i, e_i' \in \textnormal{End} \left(M(\lambda)\right)$ are given as follows:
$f_i \cdot um_\lambda=\Theta(f_i)(um_\lambda)=\varphi_\lambda f_i(u)=f_ium_\lambda$ and
$e_i' \cdot um_\lambda=\Theta(e_i')(um_\lambda)=\varphi_\lambda e_i'(u)=e_i'(u)m_\lambda$ for $u \in U_q^-(\mathfrak{g})$.
Using the same symbols for $f_i$ and $e_i'$ in both cases should create no confusion.

\begin{rem}\label{r:b} $\;$
The following conclusions are evident.
\begin{enumerit}
\item[(1)] $e_i' \cdot m_\lambda=0$ for each $i$.
\item[(2)] $f_i, e_i' \in \textnormal{End} \left(M(\lambda)\right)$ do not depend on the choice of a highest weight vector.
\item[(3)] $\varphi_\lambda : U_q^-(\mathfrak{g}) \to M(\lambda)$ is also a
$\mathcal{B}_q(\mathfrak{g})$-module isomorphism.
\item[(4)] $M(\lambda) \in \mathcal{O}\left(\mathcal{B}_q(\mathfrak{g})\right)$.
\item[(5)] $M(\lambda)$ is irreducible as a $\mathcal{B}_q(\mathfrak{g})$-module.
\end{enumerit}
\end{rem}

We emphasize Remark \ref{r:b}(5) at this point. For example, in $\mathfrak{sl_2}$-case,
$M(-n-2)$ is a
submodule of $M(n) (n \in \mathbb{Z^+}) $ with respect to the $U_q$-structure, but
as $\mathcal{B}_q$-modules $M(-n-2) \cong M(n)$.

\subl{}\label{bt}
Next, we study the $\mathcal{B}_q$-structure on $T$ modules and therefore we consider $\mathfrak{g}=\mathfrak{sl_2}$. Let $n \in \mathbb{Z^+}$. We aim to extend the $\mathcal{B}_q$-action on $M(n)$ embedded in $T(n)$ to the whole of $T(n)$.

Let $z$ be a $U_q$-generator of $T(n)$ of weight $-n-2$ and let $v$ be a highest weight vector of the Verma submodule of $T(n)$ with highest weight $n$.
Consider the $U_q^-$-decomposition
\begin{equation}\label{b9}
 T(n)=U_q^- \cdot v \oplus U_q^- \cdot z.
\end{equation}

We define $f, e' \in \textnormal{End} \left(T(n)\right)$ using (\ref{b9}) and $f, e' \in
\textnormal{End} (U_q^-)$ as follows. For $u \in U_q^-$, set
\begin{equation}\label{b10}
 \begin{array}{ll}
  f \cdot uv=fuv      &f \cdot uz=fuz\\
  e' \cdot uv=e'(u)v  &e' \cdot uz=e'(u)z.
 \end{array}
\end{equation}

The next assertion is easily seen.
\begin{lem}
Equations \eqref{b10} define a $\mathcal{B}_q$-module structure on $T(n)$.
\end{lem}

\begin{rem} $\;$

\begin{enumerit}
\item[(1)] Ker $e' =\mathbb{Q}(q)v\oplus\mathbb{Q}(q)z$
(cf. Remark \ref{r:e}).
\item[(2)] $T(n) \in \mathcal{O}(\mathcal{B}_q)$.
\end{enumerit}
\end{rem}

We call the above {\it standard $\mathcal{B}_q$-structures} on Verma and T modules.

\subl{} \label{bu}
We now consider the category $\mathcal{I}$ and collect together desirable
properties that {\it any} $\mathcal{B}_q$-structure must have with respect to the
$U_q$-structure so that a synchronization of the two would be plausible.

Let $M$ be in $\mathcal{I}$ and denote by $M^c$ the generalized
$c$-eigenspace of the quantum Casimir element $C$ on $M$. Since $M$
is $U_q^+$-finite and $C$ is central, then $M=\oplus_{r \in
\mathbb{Z}} M^{c_r}$ where
$\displaystyle{c_r=\frac{q^{r+1}+q^{-r-1}} {(q-q^{-1})^2}}$.
Clearly, $c_r=c_{-r-2}$.

\begin{defn}\label{d:i}
Let $\widetilde{\mathcal{I}}$ be the category with objects all finitely generated $U_q$-modules $M$ in the category $\mathcal{I}$ which are also equipped with a $\mathcal{B}_q$-module structure satisfying the following conditions:
\begin{enumerit}
\item [(a)] $f \in \mathcal{B}_q$ acts the same as $f \in U_q$;
\item [(b)]If $e \cdot m=0$ for a weight vector $m \in M \smallsetminus f \cdot M$ of weight $n \in \mathbb{Z}$, then $e' \cdot m=0$;
\item [(c)]If $e \cdot m \ne 0$ for a weight vector $m \in M^{c_n} \smallsetminus f \cdot M$ of weight $-n-2$ for $n \in \mathbb{Z^+}$, then there exist a
$\mathbb{Q}(q)$-subspace $T$ of $M$ such that

\begin{enumerit}
\item[(1)] $T$ is both $U_q$ and $\mathcal{B}_q$-submodule of $M$,
\item[(2)] $T \cong T(n)$ both as $U_q$ and $\mathcal{B}_q$-module where $T(n)$ is
endowed with standard $\mathcal{B}_q$-structure,
\item[(3)] $e \cdot m=e \cdot \w{m}$ for some $U_q$-generator $\w{m}$ of $T$ of weight $-n-2$.
\end{enumerit}
\end{enumerit}
Morphisms are defined to be $\mathbb{Q}(q)$-linear maps that are both $U_q$ and $\mathcal{B}_q$-morphisms.
\end{defn}

The following lemma and proposition are immediate.
\begin{lem}
Let $M$ be in $\widetilde{\mathcal{I}}, \ m \in M \smallsetminus f \cdot M,$ and $p \ge 1.$ Then
$$e'f^pm=q^{-2p}f^pe'm+\frac{1-q^{-2p}}{1-q^{-2}} f^{p-1}m. $$
\end{lem}

By the above lemma, the action of $e'$ on $M$ is determined by its action on $M \smallsetminus f \cdot M$.

\begin{prop}
\begin{enumerit}
\item The modules $\ M(r), r \in \mathbb{Z}$,
and $\ T(n), n \in \mathbb{Z^+}$, with the standard $\mathcal{B}_q$-structures belong to $\widetilde{\mathcal{I}}$.

\item Finite direct sums of modules in $\w{\mathcal{I}}$ are also in $\widetilde{\mathcal{I}}$.

\end{enumerit}
\end{prop}

\subl{}
We next aim to show:
\begin{thm}\label{t:UB}
Let $M$ be in $\widetilde{\mathcal{I}}$. Then $M=M^{(1)} \oplus \dots \oplus M^{(s)}$
for some $s>0$, where:
\begin{enumerit}
\item $M^{(j)}$ is a $U_q$ and $\mathcal{B}_q$-submodule of $M$;
\item $M^{(j)}$ is $U_q$-isomorphic either to $M(r), r \in
\mathbb{Z}$, or $T(n), \ n \in \mathbb{Z^+}$;
\item $M^{(j)}$ is $\mathcal{B}_q$-isomorphic either to $M(r), r \in
\mathbb{Z}$, or $T(n), \ n \in \mathbb{Z^+}$, with standard $\mathcal{B}_q$-structures.
\end{enumerit}
\end{thm}
\begin{proof}
If follows from Theorem \ref{t:e} and finite generation of $M$ that
$M=N^{(1)} \oplus \dots \oplus N^{(p)} \oplus \dots \oplus N^{(s)}$ where, as $U_q$-modules,
$N^{(j)} \cong M(r_j)$ for some $r_j \in \mathbb{Z}$ for
$1 \le j \le p$, and $N^{(j)} \cong T(n_j)$ for some $n_j \in \mathbb{Z}^+$ for $p+1 \le j \le s$.

For $1 \le j \le p$, let $m_j$ be a highest weight vector of $N^{(j)}$. It is clear that
$m_j \notin f \cdot M$ and $e \cdot m_j=0$. So by Definition \ref{d:i}(b), $e' \cdot m_j=0$.
It then follows from Lemma \ref{bu} that $N^{(j)}$ is stable under $e'$, i.e., $N^{(j)}$ is
a $\mathcal{B}_q$-submodule as well, and moreover it is equipped with the standard $\mathcal{B}_q$-structure.

Now consider $N^{(p+1)}$ which is $U_q$-isomorphic to $T(n_{p+1})$ and for simplicity let $n=n_{p+1}$ and $N=N^{(p+1)}$. Let $z$ be a $U_q$-generator of $N$ of weight $-n-2$ and let $v$ be a highest weight vector of the Verma submodule of $N$ with highest weight $n$.
Note that
$z \notin f \cdot M, \ e \cdot z \ne 0$, and $z \in
M^{c_n}$.
So by Definition \ref{d:i}(c), there exist a $\mathbb{Q}(q)$-subspace $T$ of $M$ which is both $U_q$ and $\mathcal{B}_q$-isomorphic to $T(n)$ with standard $\mathcal{B}_q$-action and a generator
$\w{z}$ of $T$ such that $e \cdot z=e \cdot \w{z}$. Set $w=z-\w{z}$. Then $e \cdot w=0$ and
$w$ is of weight $-n-2$.

Let $R=\oplus_{j=1, j \ne p+1}^s N^{(j)}$. Then
$M=R \oplus N$. Hence $w=w_R+w_N$ for some $w_R \in R, \ w_N \in N$, and $e \cdot w_R=0
=e \cdot w_N$. Since $w_N$ is of weight $-n-2$, it is a linear combination of $z$ and $f^{n+1}v$ and since
$e \cdot z \ne 0$, then
$w_N=\alpha f^{n+1}v$ for some $\alpha \in \mathbb{Q}(q)$. Consider
$z_0=z-w_N=z-\alpha f^{n+1}v$. Then, $z_0$ also generates $N$ and $z_0=w_R+\w{z}$. Therefore $N\subseteq R+T$ and so $M=R+T$. It is easily seen by weight consideration that the $U_q$-submodule $R \cap T$ of $T$ cannot be isomorphic neither to $T$ nor to Verma submodules of $T$ since the corresponding weight spaces of $T$ and $N$ have the same dimension and $R+T=R \oplus N$. Hence $R \cap T =\{0\}$, and so $M=R \oplus T$ as $U_q$-modules.

Therefore, we can replace $N^{(p+1)}$ in the direct sum we started with by $T$. Denote $T$ by
$T^{(p+1)}$.
Now we replace in the same way one by one the $U_q$-submodules $N^{(p+2)}, \dots,
N^{(s)}$ by $U_q$-submodules $T^{(p+2)}, \dots, T^{(s)}$ with standard
$\mathcal{B}_q$-module structure such that
$M=N^{(1)} \oplus \dots \oplus N^{(p)} \oplus T^{(p+1)} \oplus \dots \oplus T^{(s)}$.
This proves the theorem.
\end{proof}

\begin{cor}\label{c:7am}
Let $M$ be in $\w{\mathcal{I}}$. Then:
\begin{enumerit}
\item $\textnormal{Ker} \, e'=\oplus_{r \in \mathbb{Z}}
(\textnormal{Ker} \, e' \cap M_r)$;
\item $\textnormal{Ker} \, e$ has a $\qq$-basis consisting of weight vectors of the form
$f^{k}u$ where $k \in \mathbb{Z^+}$ and $u \in \, \textnormal{Ker} \, e'$.
\end{enumerit}
\end{cor}

%%%%%%%%%%%%%%%%%%%%%%%%%%%%%%%%%%%%%%%%%%%%%%%%%%%%%%%%%%%%%%%%%%%%%%

\section{Completions of Crystal Bases of Modules in the Category $\w{\mathcal{I}}$}

 Now we consider crystal bases in the sense of
Definition \ref{d:CLBCBB} for modules in $\w{\mathcal I}$ with any
$\mathcal{B}_q$-structure. Since every module in $\w{\mathcal{I}}$
has completion belonging to $\w{\mathcal{I}}$, as well as a crystal
basis, it is natural to examine the interplay between these two
concepts. We aim to define a notion of completion of crystal bases
of modules in $\w {\mathcal{I}}$ that will be compatible with the
notion of completion of modules.

\subl{}\label{lattices} The first two subsections are written for
any simple Lie algebra $\mathfrak{g}$. Let $m, n \in \mathbb{Z}^+$.
Recall (cf. \cite{K2}),
\begin{align}
&f^{(n)}_if^{(m)}_i=\q{n+m}{n}_if^{(n+m)}_i\label{i}\\
&[n]_i \in q_i^{-n+1} (1+q\ba) \ \mbox{for} \  n \ne 0\label{ii}\\
&[n]_i! \in q_i^{-\frac{n(n-1)}{2}} (1+q\ba)\label{iii}\\
&\q{m}{n}_i \in q_i^{-n(m-n)} (1+q\ba)\label{iv}
\end{align}
where $\displaystyle{\q{m}{n}_i=\frac{[m]_i[m-1]_i \dots
[m-n+1]_i}{[1]_i[2]_i \dots [n]_i}}$ for $n>0$ and
$\displaystyle{\q{m}{0}_i=1}$.

Denote by $\ba^\times$ the units in $\ba$. The next lemma is immediate and it holds for crystal lattices both
in the sense of Definition \ref{d:CLBCBB} and in the sense of Definition \ref{d:CLCB}.
\begin{lem}
Let $L$ be a crystal lattice of a module $M$ and let $\widetilde{L}=q^raL$ for some $r \in \mathbb{Z}$ and $a \in \ba^\times$. Then $\widetilde{L}$ is also a crystal lattice of $M$. If $r>0$ then $\widetilde{L} \subsetneqq L$, while if $r<0$ then $\widetilde{L} \supsetneqq L$. If $r=0$ then $\widetilde{L} = L$.
\end{lem}

\subl{} The following proposition is an analogue of a statement for
crystal bases of finite-dimensional modules, i.e., in the sense of
Definition \ref{d:CLCB}. We give a proof for reader's convenience.
\begin{prop}\label{l:iso}
Let $M$ and $M'$ be in $\mathcal{O}\left(\mathcal{B}_q(\mathfrak{g})\right)$ and let
$\varphi : M \to M'$ be a ${B}_q(\mathfrak{g})$-module isomorphism.
If $(L,B)$ is a crystal basis of $M$, then $\left(\varphi(L), \overline{\varphi}(B)\right)$
is a crystal basis of $M'$, where
$\overline{\varphi}$ is the induced map $L \slash qL \to \varphi(L) \slash q\varphi(L)$.
\end{prop}
\begin{proof}
Consider the decompositions of both $M$ and $M'$ as in (\ref{cb6}) and $\e$ and $\f$ defined as in (\ref{cb7}). It is plain that $\mbox{Ker} \, e_i' \big|_{M'} =
\varphi \left(\mbox{Ker} \, e_i' \big|_M\right)$ and both $\e$ and $\f$ commute with $\varphi$.
Thus, $\varphi(L)$
is a crystal lattice of $M'$. The restriction
$\varphi\big|_L : L \to \varphi(L)$ is an $\ba$-module
isomorphism, and so the induced map
$\overline{\varphi} : L \slash qL \to \varphi(L) \slash q\varphi(L)$ is a $\mathbb{Q}$-vector space
isomorphism. Hence $\overline{\varphi}(B)$ is a
$\mathbb{Q}$-basis of $\varphi(L) \slash q\varphi(L)$. Clearly both $\e$ and $\f$
commute with $\overline{\varphi}$.
\end{proof}

For $\lambda \in P$, let $m_\lambda$ be a highest weight vector of the Verma module $M(\lambda)$ and let $\varphi_\lambda : U_q^-(\mathfrak{g}) \to M(\lambda)$
be a $\mathcal{B}_q(\mathfrak{g})$-module isomorphism such that
$\varphi_\lambda(u)=um_\lambda$. Let $L^{(\lambda)}$ be the $\ba-\mbox{span of} \
\{\widetilde{f_{i_1}} \dots \widetilde{f_{i_p}} \cdot m_\lambda | \ 1 \le i_j \le l;
p \in \mathbb{Z}^+\}$ and
$B^{(\lambda)}=\{\widetilde{f_{i_1}} \dots \widetilde{f_{i_p}} \cdot \overline{m_\lambda}
| \ 1 \le i_j \le l; p \in \mathbb{Z}^+\}$, where
$\overline{m_\lambda}=m_\lambda+q L^{(\lambda)}$. Since $\left(L(\infty), B(\infty)\right)$ is a crystal basis
of $U_q^-(\mathfrak{g})$, it immediately follows:
\begin{cor}\label{r:a}
$(L^{(\lambda)}, B^{(\lambda)})$ is a crystal basis of the Verma module $M(\lambda)$. Every crystal basis of
$M(\lambda)$ is a nonzero scalar multiple of $(L^{(\lambda)}, B^{(\lambda)})$.
\end{cor}

\begin{ex}\label{e:V}
If $\mathfrak{g}=\mathfrak{sl_2}$ and $r \in \mathbb{Z}$, then any crystal basis of the Verma module $M(r)$ is of the form
$(L^{(r)}, B^{(r)})$ where $L^{(r)}=\oplus_{k \in \mathbb{Z^+}}
\ba f^{(k)}m_r$ and
$B^{(r)}=\{f^{(k)}\overline{m_r} | \ k \in \mathbb{Z^+}\}$ for a highest weight vector $m_r$.
\end{ex}

\subl{}\label{r:bb}

From this point on, we focus exclusively on the
$U_q(\mathfrak{sl_2})$-case.

We consider $T$ modules. Recall that
$\left(L(\infty)^{\oplus 2}, B(\infty)^{\oplus 2}\right)$ is a crystal basis of
$U_q^- \oplus U_q^-$, where $L(\infty)^{\oplus 2}=
L(\infty) \oplus L(\infty)$, $B(\infty)^{\oplus 2}=\left(B(\infty) \times 0\right)
\cup \left(0 \times B(\infty)\right) \subset
L(\infty)^{\oplus 2} \slash qL(\infty)^{\oplus 2}$, $L(\infty)=\oplus_{k \in \mathbb{Z^+}} \ba f^{(k)}
\cdot 1$ and
$B(\infty)=\{f^{(k)} \cdot \overline{1} | \ k \in \mathbb{Z^+}\}$ where $\overline{1}=1+qL(\infty)$.
The next corollary follows immediately from considering the map $\psi_n :
U_q^- \oplus U_q^- \to T(n)$
given by $(u_1,u_2) \mapsto u_1v+u_2z$ which is a $\mathcal{B}_q$-module
isomorphism where $U_q^- \oplus U_q^-$ is equipped with the obvious
$\mathcal{B}_q$-structure.

\begin{cor}
Let $z$ be a $U_q$-generator of $T(n)$ of weight $-n-2$ and let $v$ be a highest weight vector of Verma submodule of $T(n)$ with highest weight $n$. Set
\[
 L^{\{n\}}=\left(\oplus_{k \in \mathbb{Z^+}} \ba f^{(k)}v\right) \oplus \left(\oplus_{k \in \mathbb{Z^+}}
\ba f^{(k)}z\right) \qquad and
\]
\[
 B^{\{n\}}=\{f^{(k)}v \, \textnormal{mod} \, q L^{\{n\}} |
k \in \mathbb{Z^+}\} \cup \{f^{(k)}z \, \textnormal{mod} \, q L^{\{n\}} | \ k \in \mathbb{Z^+}\}.
\]
Then $(L^{\{n\}}, B^{\{n\}})$ is a crystal basis of $T(n)$. Moreover, if $(L,B)$ is any crystal basis of $T(n)$, then there exist $z$ and $v$ as above such that $(L,B)$ is of the form $(L^{\{n\}}, B^{\{n\}})$.
\end{cor}

The following theorem now easily follows from Theorem
\ref{t:UB}, Corollary \ref{r:a} and Corollary \ref{r:bb}.

\begin{thm}\label{t:cor}
Let $M$ be a module in the category $\w{\mathcal{I}}$. Then $M$
belongs to the category $\mathcal{O}(\mathcal{B}_q)$ and there exist unique integers $r_1,...,r_m$ and unique nonnegative integers $n_1,...,n_t$ such that $M$ has a crystal basis $(L,B)$ which is a direct sum of crystal bases of the form $(L^{(r_i)}, B^{(r_i)})$, for $i=1,...,m$, and crystal bases
of the form $(L^{\{n_j\}}, B^{\{n_j\}})$, for $j=1,...,t$. Furthermore, any
crystal basis of $M$ is isomorphic to $(L,B)$.
\end{thm}

\begin{defn}
We call the crystal basis $(L,B)$ from the previous theorem a standard crystal basis of $M$.
\end{defn}
Evidently, the integers $r_1,...,r_m$ and $n_1,...,n_t$ depend only on $M$ and not on $(L,B)$.
Every crystal basis of an indecomposable module in $\w{\mathcal{I}}$ is obviously a standard crystal basis.

%%%%%%%%%%%%%%%%%%%%%%%%%%%%%%%%%%%%%%%%%%%%%%%%%%%%%%%%%%%%%%%%%%%%%%%%%%%%%%%%%

\subl{}

For $M$ in $\w{\mathcal{I}}$, recall the decomposition $M=\oplus_{k \ge 0} f^{(k)} \mbox{Ker} \, e'$ (see \ref{d:CLBCBB}). By Corollary \ref{c:7am},
$\mbox{Ker} \, e'=\oplus_{n \in \mathbb{Z}} \left(\mbox{Ker} \, e'
\cap M_n\right)$. Thus we bring the weight spaces of $M$ in the picture and have:
\begin{equation}\label{ccl1}
 M=\bigoplus_{\substack{k \ge 0\\ n \in \mathbb{Z}}} f^{(k)} \left(\mbox{Ker} \,
e' \cap M_n\right).
\end{equation}

For a crystal lattice $L$ of $M$, set $L_n=L \cap M_n$. Then $L= \oplus_{n\in \mathbb{Z}}L_n$. Also, set $L^{e} = \{m \in L \,| \, em=0\} = L \cap M^e \
\mbox{and} \ L_n^{e} = L^{e} \cap M_n = L \cap M^e_n.$

\begin{lem}\label{ccl2} Let $M$ be in $\w{\mathcal{I}}$ and let
$L$ be a crystal lattice of $M$. Then ${\widetilde{f}}^{n+1}(M_n^e)
\subseteq M_{-n-2}^e$ and $\w{f}^{n+1}(L_{n}^e) \subseteq L_{-n-2}^e$
for $n \in \mathbb{Z}^+$.
\end{lem}
\begin{proof}
Let $m \in M_n^e$. By Corollary \ref{c:7am}, $m=\sum_{k \ge 0} f^{(k)}u_k$
where $f^{(k)}u_k \in M^e$ and $u_k \in \, \mbox{Ker} \, e' \cap M_{n+2k}$.
Thus ${\widetilde{f}}^{n+1} m=\sum_{k \ge 0} {\widetilde{f}}^{n+1} f^{(k)} u_k=
\sum_{k \ge 0} f^{(n+1+k)} u_k \in M_{-n-2}$. Furthermore, utilizing (\ref{i}) and commutation relations,
$$e\w{f}^{n+1}f^{(k)}u_k
\displaystyle{
=\frac{1}{\q{n+k+1}{n+1}}ef^{(n+1)}f^{(k)}u_k
=\frac{1}{\q{n+k+1}{n+1}}\left(f^{(n+1)}e+f^{(n)}\frac{q^{-n}t-q^{n}t^{-1}}{q-q^{-1}}\right)f^{(k)}u_k
=0}$$ for all $k \ge 0$.
Therefore $\w{f}^{n+1}m \in M^e$, and so
${\widetilde{f}}^{n+1}(M_n^e) \subseteq M_{-n-2}^e.$
In addition, since $\widetilde{f}(L) \subseteq L$, then
${\widetilde{f}}^{n+1}(L_{n}^e) \subseteq L_{-n-2}^e$.
\end{proof}

\begin{defn}\label{d:CCL}
Let $L$ be a crystal lattice of a module $M$ in $\w{\mathcal{I}}$. We say $L$ is complete if  $\
{\widetilde{f}}^{n+1} : L_{n}^e \to L_{-n-2}^e$ is bijective for all $n \in \mathbb{Z}^+.$
A crystal lattice $\widetilde{L}$ of the completion $C(M)$ of $M$ is said to be
a completion of a crystal lattice $L$ of $M$ provided
\begin{enumerit}
\item $L^e \subseteq \w{L}$ and $\w{L} \ns L$;
\item $\w{L} \cap M = \w{L} \cap L$;
\item $\w{L} \slash (\w{L} \cap L)$ is a crystal lattice of $C(M) \slash M$
(to be precise, its image under $\w{L} \slash (\w{L} \cap L) \hookrightarrow \break C(M) \slash M$).
\end{enumerit}
\end{defn}

\begin{rem}
\begin{enumerit}
\item Unlike the case of modules where $C(M)$ contains $M$, one cannot expect $\w{L}$ to contain
$L$. If we consider crystal lattices of a reducible Verma module and its irreducible Verma $U_q$-submodule, denoted by $L_1$ and $L_2$ respectively, although $L_1$ is isomorphic to $L_2$, when talking about completions we are
interested in their exact relationship and indeed $L_2 \nsubseteq L_1$. Moreover, we observe
that there is no $p \in \mathbb{Z}$ such that $q^pL_2 \subseteq L_1$. In fact, the actions of the Kashiwara operators $\w{e}$ and $\w{f}$ on $L_1$ and $L_2$ are different, i.e., $\w{e}_M \neq (\w{e}_{C(M)})|_M$ and $\w{f}_M \neq (\w{f}_{C(M)})|_M$.
\item Condition (iii) in Definition \ref{d:CCL} is a very natural connection to
expect between crystal lattices and completions. Also, we note that
$C(M) \slash M$ is a finite-dimensional $U_q$-module, thus condition
(iii) gives a connection between crystal lattices arising from the
$\mathcal{B}_q$-structures, i.e., in the sense of Definition
\ref{d:CLBCBB}, with crystal lattices arising from the
$U_q$-structures, i.e., in the sense of Definition \ref{d:CLCB}.
\item The condition $\w{L} \ns L$, i.e., $\w{L}$ is not properly contained in $L$, follows from condition (iii) in the case that $L$ is a crystal lattice of a module $M$ that is not complete.
\end{enumerit}
\end{rem}

\begin{prop}
Let $L$ be a crystal lattice of a module $M$ in $\w{\mathcal{I}}$ and let $\w{L}$ be a completion of $L$. Then $\w{L}^e_{-n-2}=L^e_{-n-2}$ for all $n \in \mathbb{Z^+}$.
\end{prop}

\begin{proof}
Let $n \in \mathbb{Z^+}$. We note that since $C(M)/M$ is a
finite-dimensional $U_q$-module, then
$\left(C(M)/M\right)^e_{-n-2}=0,$ and so $C(M)^e_{-n-2}=M^e_{-n-2}$.
Now, by Definition \ref{d:CCL}(ii), $\w{L}^e_{-n-2}=\w{L} \cap
C(M)^e_{-n-2}=\w{L} \cap M^e_{-n-2}= (\w{L} \cap M)\cap M^e_{-n-2} =
(\w{L} \cap L) \cap M^e_{-n-2}=\w{L}^e_{-n-2}\cap L^e_{-n-2}$.
Hence, $\w{L}^e_{-n-2}\subseteq L^e_{-n-2}$. The claim now follows
from Definition \ref{d:CCL}(i).
\end{proof}

\subl{} Let $(L,B)$ be a crystal basis of $M$. Set $B_n=L_n/qL_n
\cap B \,$ where $\,L_n/qL_n=\{m \, \text{mod} \, qL \,|\, m \in
L_n\}$. Then $B= \sqcup_{n\in \mathbb{Z}}B_n$. Also set
$B_n^e=(L_n/qL_n)^e \cap B \,$ where $\,(L_n/qL_n)^e=\{m \,
\text{mod} \, qL \,|\, m \in L_n^e\}$.
\begin{lem}
Let $M$ be in $\w{\mathcal{I}}$ and let
$(L,B)$ be a crystal basis of $M$. Then:
\begin{enumerit}
\item ${\widetilde{f}}^{n+1}(B_n^e)\subseteq B_{-n-2}^e$ for $n \in \mathbb{Z}^+$;
\item If $L$ is complete, then $\widetilde{f}^{n+1} : B_{n}^e \to B_{-n-2}^e$ is bijective for all $n \in \mathbb{Z}^+.$
\end{enumerit}
\end{lem}

\begin{proof}
Part (i) of the Lemma follows from Lemma \ref{ccl2} and from $B$ being invariant under $\w{f}$. It is easily seen from Corollary \ref{c:7am} that $B_{n}^e$ is a basis of $(L_n/qL_n)^e$, thus implying part (ii).
\end{proof}

Let $\displaystyle{\overline\varphi : \, \w L/q\w L \to\frac{\w L/(L
\cap \w L)} {q\left(\w L/(L \cap \w L)\right)}}$ be the induced
$\mathbb{Q}$-epimorphism from the canonical $\mathbb{A}$-epimorphism
$\varphi: \, \w{L} \to \w{L}/{(L \cap \w{L})}$.

\begin{defn}\label{d:CCB}
A crystal basis $(L,B)$ of $M$ in $\w{\mathcal{I}}$ is said to be complete if the crystal
lattice $L$ is complete.

\noindent A crystal basis $(\widetilde{L},\widetilde{B})$ of $C(M)$
is a completion of a crystal basis $(L,B)$ if
\begin{enumerate}
\item $\widetilde{L}$ is a completion of $L$,
\item $\left(\w L/(L \cap \w L), \ \overline\varphi(\w B) \right)$
is a crystal basis of $C(M)/M$,
\end{enumerate}
\end{defn}

The next theorem verifies that the above definition of
completeness for crystal bases is compatible with the notion of
completeness for modules.

\begin{thm}\label{t:CLC}
Let $(L,B)$ be a crystal basis of a module $M$ in $\w{\mathcal{I}}$.
Then $(L,B)$ is complete if and only if $M$ is complete.
\end{thm}
\begin{proof}
Let $M$ be complete and let $n \in \mathbb{Z}^+$. Then $f^{n+1} :
M_n^e \to M_{-n-2}^e$ is bijective. Since $M$ is in
$\w{\mathcal{I}}$, its weight spaces are finite-dimensional. Hence,
\begin{equation}\label{ccl3}
 \dim_{\mathbb{Q}(q)} M_n^e = \dim_{\mathbb{Q}(q)} M_{-n-2}^e.
\end{equation}
$L$ is a crystal lattice of $M$, so there is an $\ba$-basis of $L$ that is also a
$\mathbb{Q}(q)$-basis of $M$. This basis is made up of weight vectors
of the form $f^{(k)}u$ where $k \ge 0$ and $u \in \, \mbox{Ker} \, e'$. If
a nonzero sum of such vectors is annihilated by $e$, then also each of them is
(see Theorem \ref{t:UB} and Corollary \ref{c:7am}).
Thus, a $\mathbb{Q}(q)$-basis of
$M_n^e$ is an $\ba$-basis of $L_n^e$.
By (\ref{ccl3}),
$\mbox{rank}_{\ba} \, L_{n}^e = \mbox{rank}_{\ba} \, L_{-n-2}^e$. Now, since $\ba$ has the invariant dimension property and since $\w{e} \w{f}=1$ on $M$,
then ${\widetilde{f}}^{n+1} : L_{n}^e \to L_{-n-2}^e$ is bijective.
The converse is true since $f$ is injective on $M$, and the above steps can be reversed for a fixed $n$.
\end{proof}

\begin{ex}\label{e:ccl1}
The modules $T(n)$ for $n \in \mathbb{Z^+}$ and $M(n)$ for $n \ge
-1$ are complete, so their crystal bases defined in \ref{r:bb} and
\ref{r:a}, respectively, are complete by the theorem. This can also
be seen directly from the Definitions \ref{d:CCL} and \ref{d:CCB}.
\end{ex}

\begin{cor}
A complete crystal basis $(L,B)$ of a module $M$ in $\w{\mathcal{I}}$ is a
completion of itself. Furthermore, $L$ is the only completion of itself.
\end{cor}

\begin{proof}
By the previous theorem, $(L,B)$ being complete implies $M=C(M)$. It is obvious that a crystal basis
of a complete module satisfies Definitions \ref{d:CCL} and \ref{d:CCB} for being a completion of itself.
Now assume $\w{L}$ is another completion of $L$. Then $\w{L}=\w{L} \cap C(M)=\w{L}\cap M=\w{L}\cap L$,
thus $\w{L}\subseteq L$. The desired conclusion follows from the condition $\w{L} \ns L$.
\end{proof}

\subl{}
The rest of the section is devoted to the proof of the following theorem.
\begin{thm}\label{t:exuni}
A standard crystal basis $(L,B)$ of any module $M$ in $\w{\mathcal{I}}$ has a
completion. Moreover, there exists a unique standard crystal lattice which is a completion of $L$.
\end{thm}

By Theorem \ref{t:cor},
every module $M$ in category $\w{\mathcal{I}}$ has a standard crystal basis, thus it has a crystal basis that can be completed as in the above theorem.

We first prove:

\begin{prop}\label{p:ind}
A crystal basis $(L,B)$ of any indecomposable module $M$ in $\w{\mathcal{I}}$
has a completion. Furthermore, $L$ has a unique completion.
\end{prop}

\begin{proof}
By Theorem \ref{t:CLC} and Corollary \ref{t:CLC}, in order to show existence, we need to
consider only Verma modules in $\w{\mathcal{I}}$ which are not
complete. Let $M$ be a Verma module with highest weight $-n-2$ for
some $n \in \mathbb{Z^+}$ and let $L$ be its crystal lattice. Then
$L=\oplus_{k\ge0}\ba f^{(k)}m_0$ for some highest weight vector
$m_0$ of $M$ and $L$ is clearly not complete. Since $M$ has a
completion $C(M)$ and $f^{n+1} : C(M)_{n}^e \to
C(M)_{-n-2}^e=M_{-n-2}^e$ is bijective, there exists $\w m$ in
$C(M)_n^e$ such that $f^{(n+1)}\w m=m_0$. Let $\w L=\oplus_{k\ge
0}\ba f^{(k)}\w m$ and $\w{B}=\{f^{(k)}\w{m}+q\w{L} | \, k \in
\mathbb{Z}^+\}$. Then $(\w L, \w B)$ is a crystal basis of $C(M)$
and we claim it is a completion of $(L,B)$.

It follows from subsection \ref{lattices} that
$$L=\oplus_{k \ge 0} \ba f^{(k)}f^{(n+1)} \w m
   =\oplus_{k \ge 0} \ba \q{n+1+k}{k} f^{(n+1+k)}\w m
   =\oplus_{k \ge 0} \ba q^{-k(n+1)} f^{(n+1+k)}\w m.$$
The last equality is a consequence of the elements of $1+q\ba$ being units in $\ba$.

We have $L^e=\ba m_0 \subseteq \w{L}$,
$\ \w L \cap L = \oplus_{k\ge 0} \ba f^{(n+1+k)}\w m$, and $\w L \slash (\w L \cap L) = \oplus_{k=1}^n \ba f^{(k)}\overline{m}$ where $\overline{m}= \w m+q(\w L \cap L)$. Thus the conditions (i)-(iii) of Definition \ref{d:CCL} are satisfied. It is now easy to check condition (2) of Definition \ref{d:CCB}.

Uniqueness follows from parts (i) and (ii) of Definition \ref{d:CCL}.
\end{proof}

\begin{rem}\label{v:down}
\begin{enumerit}
\item Considering the lattice $\widetilde L $ from the proof of previous proposition, notice that $\widetilde L \cap M= \oplus_{k \ge 0}\ba f^{(n+1+k)}\w m$ which is a
proper $\ba$-submodule of the crystal lattice $L=\oplus_{k \ge 0} \ba q^{-k(n+1)}f^{(n+1+k)}\w m$ of $M$, i.e., in order to complete $L$ to a crystal lattice of $C(M)$ that is invariant under Kashiwara operators $\widetilde e_{C(M)}$ and $\widetilde f_{C(M)}$, $\ L$ has to be made ``thinner'' as an $\ba$-lattice in a particular way.

\item Using the proof of the previous proposition, we can also see the following.
Let $n \in \mathbb{Z}^+$ and let $M=M(n)$ be the Verma module
with highest weight $n$ and $M'$ its
unique Verma submodule. Let $L$ be any crystal lattice of $M$ and $L'$ any crystal lattice of $M'$. Then,
there exists $r \in \mathbb{Z}^+$ such that $L \slash (L \cap q^{-r}L')$ is isomorphic to
a crystal lattice of the irreducible $U_q$-module $M \slash M'$.
\end{enumerit}
\end{rem}

\subl{}

\begin{prop}\label{p:sum}
Direct sum of completions of crystal bases is a completion of direct
sum of crystal bases.
\end{prop}
\begin{proof}
For $i=1,2$, let $(L_i,B_i)$ be a crystal basis of $M_i$ and $(\w{L_i},\w{B_i})$ its
completion. We claim that $(\w{L_1},\w{B_1}) \oplus (\w{L_2},\w{B_2})$ is a completion
of $(L_1,B_1) \oplus (L_2,B_2)$. Since both crystal bases and completions of
modules in $\w{\mathcal{I}}$ respect direct sums, $(\w{L_1},\w{B_1}) \oplus (\w{L_2},\w{B_2})=(\w{L_1} \oplus
\w{L_2},\w{B_1} \sqcup \w{B_2})$ is a crystal basis of $C(M_1) \oplus C(M_2) = C(M_1
\oplus M_2)$. Condition (i) of Definition \ref{d:CCL} is evident
since $(L_1 \oplus L_2)^e=L_1^e \oplus L_2^e \subseteq \w{L_1}
\oplus \w{L_2}$. Furthermore, $\w{L_1} \oplus \w{L_2} \ns L_1 \oplus
L_2$. By definition, $\w{L_i} \cap M_i = \w{L_i} \cap L_i$. Clearly,
for $i \ne j$, $M_i \cap \w{L_j} =0$. Therefore, $(\w{L_1} \oplus
\w{L_2}) \cap (M_1 \oplus M_2) = (\w{L_1} \oplus \w{L_2}) \cap (L_1
\oplus L_2)$, so condition (ii) is satisfied, as well. Next,
$$(\w{L_1} \oplus \w{L_2}) \slash \big((\w{L_1} \oplus \w{L_2}) \cap (L_1 \oplus L_2)\big)=\break
 (\w{L_1} \oplus \w{L_2}) \slash \big((\w{L_1} \cap L_1) \oplus (\w{L_2} \cap L_2)\big)
 \cong \w{L_1} \slash (\w{L_1} \cap L_1) \ \oplus \ \w{L_2} \slash (\w{L_2} \cap L_2).$$
On the other hand,
$$C(M_1 \oplus M_2) \slash (M_1 \oplus M_2)
= \big(C(M_1) \oplus C(M_2)\big) \slash \big(M_1 \oplus M_2\big)
\cong C(M_1) \slash M_1 \ \oplus \ C(M_2) \slash M_2.$$
Condition (iii) now follows from crystal lattices respecting direct sums. Furthermore, it is easy to see that $\bar{\varphi}(\w{B_1} \sqcup \w{B_2}) = \bar{\varphi}(\w{B_1}) \sqcup \bar{\varphi}(\w{B_2})$ proving condition (2) of Definition \ref{d:CCB}.
\end{proof}

\begin{lem}
Let $M$ be a module in $\w{\mathcal{I}}$ and suppose $L_i, L_i', i=1,\dots, k$, are crystal lattices of indecomposable submodules of $M$ such that $L_1\oplus L_2 \oplus \cdots \oplus L_k = L'_1\oplus L'_2\oplus \cdots\oplus L_k'$. If $\w{L_i}, \w{L_i'}$ are the completions of $L_i, L_i'$, respectively, then  $\w{L_1}\oplus \w{ L_2} \oplus \cdots \oplus \w{L_k} = \w{L'_1}\oplus \w{L'_2}\oplus \cdots\oplus \w{L_k'}$.
\end{lem}
\begin{proof}
It suffices to consider the case when $M$ is a generalized eigenspace of the Casimir element with eigenvalue $c_n$. If all of the $L_i$ are complete, then so is $L_1\oplus L_2 \oplus \cdots \oplus L_k$ and we are done using Proposition \ref{p:sum} and Corollary \ref{e:ccl1}. To simplify the notation we consider only the case when $k=2$, $L_1$ is a crystal lattice of a Verma submodule of $M$ with highest weight $n$, and  $L_2$ is a crystal lattice of a Verma submodule of $M$ of highest weight $-n-2$. The proof of the general case is similar, but with more involved notation. It is clear that either $L_1'$ or $L_2'$ is a crystal lattice of a  Verma submodule of $M$ with highest weight $n$, so we assume it is $L_1'$. Furthermore, it is then obvious from the condition $L_1\oplus L_2=L'_1\oplus L'_2$ that we must have $L_1 = L_1' = \w{L_1} = \w{L_1'}$. Let $v\in {\rm Ker}\,  e'\cap M_n$ and $u,u'\in {\rm Ker}\,  e'\cap M_{-n-2}$ be such that $L_1, L_2, L_2'$ are the $\ba$-spans of $f^{(i)}v, f^{(i)}u, f^{(i)}u' (i \ge 0)$, respectively. Clearly $u'=af^{(n+1)}v+bu$ where $a,b\in \ba$. Since $u \in L'_1\oplus L'_2$, then $b \in \ba^\times$. Let $m,m'$ be the unique elements in $C(M)$ such that $f^{(n+1)}m=u$ and $f^{(n+1)}m'=u'$.
Then $\w{L_2}= \sum_{i\ge0}\ba f^{(i)}m$ and $\w{L'_2}=\sum_{i\ge 0}\ba f^{(i)}m'$ by the proof of Proposition \ref{p:ind}. Also, $f^{(n+1)}m'=f^{(n+1)}(av+bm)$. Since $M$ is $U_q^-$-torsion free, then $m'=av+bm \in \w{L_1}\oplus \w{L_2}$ and $m=b^{-1}(m'-av) \in \w{L_2'}\oplus \w{L_1'}$ proving the lemma in this case.
\end{proof}

\subl{}
It is convenient to have the following definition.
\begin{defn}
Let $M_1$, $M_2$ be modules in $\w{\mathcal{I}}$ and let $(L_1,B_1)$, $(L_2,B_2)$ be crystal bases of $M_1, M_2$, respectively. We say that $(L_1,B_1)$ is strongly isomorphic to $(L_2,B_2)$ if there
exists a $B_q$-isomorphism $\varphi: M_1 \to M_2$ such that
\begin{enumerit}
\item  $\varphi$ induces an isomorphism of the crystal bases $(L_1,B_1)$ and $(L_2,B_2)$;
\item $\varphi$ is weight-preserving;
\item $\varphi(\mbox{Ker} \, e|_{M_1})=\mbox{Ker} \, e|_{M_2}.$
\end{enumerit}
\end{defn}

For example, any two standard crystal basis of a module in $\w{\mathcal{I}}$ are strongly isomorphic.

\noindent {\it Proof of Theorem \ref{t:exuni}}.
By Definition \ref{t:cor}, a standard crystal basis is a direct sum of crystal bases of the form $(L^{(r_i)}, B^{(r_i)})$, for $i=1,...,m$, and crystal bases of the form $(L^{\{n_j\}}, B^{\{n_j\}})$, for $j=1,...,t$,
where the integers $r_i$ and $n_j$ depend on $M$, i.e., these are the highest weights of Verma and $T$ modules present in a decomposition of $M$. The first claim of the theorem now follows from Propositions \ref{p:ind} and \ref{p:sum}.

We now prove the uniqueness of completion of $L$. Since $L$ is a standard crystal lattice of $M$, then $L=L_1\oplus \cdots \oplus L_k$ for some crystal lattices $L_i, i=1,...,k,$ of indecomposable $U_q$-submodules $M_i$ of $M$ such that $M=M_1\oplus \cdots \oplus M_k$. If $\w{L_i}$ denotes the completion of $L_i$, then $\w{L_1}\oplus \cdots \oplus \w{L_k}$ is a completion of $L$. Let $\w{L'}$ denote another standard completion of $L$. Then $\w{L'}$ is a crystal lattice of $C(M)=C(M_1)\oplus \cdots \oplus C(M_k)$, thus it is a direct sum of standard crystal lattices $L_i'', i=1,...,k,$ of indecomposable submodules $M_i''$ of $C(M)$ such that $C(M)=M_1'' \oplus \cdots \oplus M_k''$ and such that $L_i''$ is strongly isomorphic to $\w{L_i}, i=1,...,k$ (up to reordering). Since $L, L_1'',..., L_k''$ are standard crystal lattices and $L_1''\oplus \cdots \oplus L_k''$ is a completion of $L$, it is easily seen that there exist crystal lattices $L'_i, i=1,...,k,$ strongly isomorphic to $L_i, i=1,...,k$, respectively, such that $L_i''$ is a completion of $L'_i$ and $L=L'_1 \oplus \cdots \oplus L'_k$. The conclusion now follows from the previous lemma.
\hfill $\Box$

%%%%%%%%%%%%%%%%%%%%%%%%%%%%%%%%%%%%%%%%%%%%%%%%%%%%%%%%%%%%%%%%%%%%%%%%%%%%%%%%%%%%%%

\section{Constructions of Completions of Crystal Lattices}

We now construct a completion of a crystal lattice modifying Deodhar's construction of completion
for modules in $\w{\mathcal{I}}$. By Theorems \ref{t:UB} and \ref{t:CLC} and
Proposition \ref{p:sum}, it suffices to consider crystal bases of Verma modules $M(-n-2)$ for $n \in \mathbb{Z^+}$.

\subl{} \label{v:up}

We start by recalling Deodhar's construction for
$U_q(\mathfrak{sl_2})$ (cf. \cite{D}, \cite{Z}). Denote by
$\mathcal{A}\left(U_q\right)$ the category of $U_q$-weight modules
that are $U_q^-$-torsion free. For $M$ in
$\mathcal{A}\left(U_q\right)$, let $S_M$ be the set of formal
symbols $S_M=\{f^{-r}m| \ r \in \mathbb{Z}^+, m \in M\}$. Define an
equivalence relation $\sim$ on $S_M$ by $f^{-r} m \sim f^{-k} m'$
iff $f^k m = f^r m'$. Set $D(M)=S_M \slash \sim$. Then $D(M)$ has a
vector space structure. For $z \in U_q$ and $r \in \mathbb{Z}^+$,
there exist $u \in U_q$ and $s \in \mathbb{Z}^+$ such that $f^s
z=uf^r$. Now, for $f^{-r} m \in D(M)$ and $z \in U_q$, define $z
\cdot f^{-r} m=f^{-s} um$. This action makes $D(M)$ a $U_q$-module.

Any nonzero element $v \in D(M)$ can be uniquely expressed as $f^{-n} m
\ (n \in \mathbb{Z}^+, m \in M)$ where $n$ is minimal with respect to this property. This expression
is called the minimal expression for $v$. Clearly, $f^{-n} m$ is a minimal expression iff
$m \notin f M$. Also, $v \in M$ iff $n=0$ in the minimal expression for $v$.

Set $C(M)=\{v \in D(M)| \ e^{k} v=0 \ \mbox{for some} \ k>0 \}$.

\begin{thm} \label{t:d}(cf. \cite{D}, \cite{Z})
Let $M$ be a module in the category $\mathcal{I}$.
\begin{enumerit}
\item $C(M)$ is a completion of $M$.

\item Let $v \in D(M) \smallsetminus M$ be a weight vector such that $t v=q^a v$
for some $a \in \mathbb{Z}$ and let $v=f^{-n} m$ be the minimal expression for $v$. Then
$v \in C(M)$ iff the following two conditions are satisfied:
(a) \ $a=n-j$ for some $j>0$, and
(b) $e^k f^{k-j} m=0$ for some $k \ge j$.
\end{enumerit}
\end{thm}

The next corollary gives basis vectors for
the completion $C(M)$ of an irreducible Verma module $M$.

\begin{cor}
Let $M$ be the Verma module with highest weight $-n-2$ for
$n \in \mathbb{Z^+}$ and let $m_0$ be a highest weight vector of $M$.
Then the completion of $M$ is $\displaystyle{C(M)=\oplus_{k \ge -n-1} \mathbb{Q}(q)f^km_0}$.
\end{cor}

\begin{proof}
Applying Deodhar's construction on $M=\oplus_{k \ge 0} \mathbb{Q}(q)f^km_0$, we obtain
$D(M)=\oplus_{k \in \mathbb{Z}} \mathbb{Q}(q)f^km_0$ because $f$ is injective and
$f^{-r}f^km_0=f^{k-r}m_0$ for $r, k \in \mathbb{Z}^+$.
Since $M \subset C(M)$, we need to consider only $f^{-k}m_0$ for $k>0$.
By Theorem \ref{t:d}(ii), $f^{-k}m_0 \in C(M)$ iff
(a)\ $-n-2+2k=k-j$ for some $j>0$,
and (b) $e^pf^{p-j}m_0=0$ for some $p \ge j$. This is equivalent to $k<n+2$ ($p=j=n+2-k$).
\end{proof}

\subl{} \label{fac}

It is convenient to introduce the following formal notation.
Set $[-k]! = [-1] \dots [-k]$ and $\displaystyle{f^{(-k)}m=\frac{f^{-k}m}{[-k]!}}$ for $k>0$ and $m \in M$.
Then $[-k]!=(-1)^k[k]!$ and $\displaystyle{f^{(k)}m=\frac{f^km}{[k]!}}$ for $k \in \mathbb{Z}$ and  $m \in M.$

For $i, k \in \mathbb{Z}^+$ and $m \in M$,
\begin{equation}\label{cc13}
 f^{(i)}f^{(-k)}m=\frac{f^i}{[i]!} \frac{f^{-k}m}{[-k]!}=\frac{f^{i-k}m}{[i]![-k]!}
                 =\frac{[i-k]!}{[i]![-k]!} f^{(i-k)}m.
\end{equation}
Thus, $$C(M)=\oplus_{k \ge -n-1} \mathbb{Q}(q)f^{(k)}m_0=\oplus_{i \ge 0}
\mathbb{Q}(q)f^{(i)}f^{(-n-1)}m_0,$$
where we follow the notation of Corollary \ref{t:d}.

\subl{}
We now construct directly from the crystal lattice of an irreducible Verma module
a crystal lattice of its completion.

\begin{lem}\label{l:*} $\:$

\begin{enumerit}
\item[(a)] If  $i \le n$,  then $\displaystyle{1=\frac{1}{1+q^{-(i-n-1)}} \, a^*}$ for some
$a^* \in \ba^\times$.
\item[(b)] If $i > n$,  then  $\displaystyle{q^{i-n-1}=\frac{1}{1+q^{-(i-n-1)}} \, a^*}$ for some
$a^* \in \ba^\times$.
\end{enumerit}
\end{lem}
\begin{proof}
(a) \ For $i \le n, \ 1+q^{-i+n+1} \in \ba^\times$.

\noindent (b) \ For $i > n$, $\displaystyle{\frac{1}{1+q^{-(i-n-1)}}=\frac{q^{i-n-1}}{q^{i-n-1}+1}}$,
so let $a^*=1+q^{i-n-1} \in \ba^\times$.
\end{proof}
\begin{prop}\label{p:neg}
Let $M$ be the Verma module with highest weight $-n-2$ for some $n
\in \mathbb{Z^+}$ and let $m_0$ be a highest weight vector of $M$.
Then
\begin{enumerit}
\item the $\ba$-module $\w L$ generated by
$\displaystyle{\left\{\frac{2q^{n(i-n-1)}}{1+q^{-(i-n-1)}} \, f^{(i-n-1)}m_0 | \ i \ge 0\right\}}$
is a crystal lattice of $C(M)$;
\item $\w L$ is a completion of the crystal lattice $L = \oplus_{k\ge0}\ba f^{(k)}m_0$ of $M$.
\end{enumerit}
\end{prop}

\begin{proof}
\begin{enumerit}
\item By Corollary \ref{t:d}, $C(M)$ is a Verma module with a highest weight vector $f^{(-n-1)}m_0$. Thus,
 a crystal lattice of $C(M)$ is of the form $\widetilde L=\oplus_{i \ge 0}\ba r(q)f^{(i)}f^{(-n-1)}m_0$ for
$r(q) \in \mathbb{Q}(q)$. Utilizing (\ref{iii}), we have the following:
\begin{align*}
f^{(i)}f^{(-n-1)}m_0&=(-1)^{n+1}\frac{[i-n-1]!}{[i]![n+1]!} \, f^{(i-n-1)}m_0\\
                    &=(-1)^{n+1}\frac{[i-n-1]!}{q^{-\frac{1}{2}[(i-1)i+n(n+1)]}(1+aq)}
                      f^{(i-n-1)}m_0 \, \ \mbox{for some} \ a \in \ba.
\end{align*}
\textit{Case 1:} \ $i \ge n+1$
\begin{align*}
 f^{(i)}f^{(-n-1)}m_0&=(-1)^{n+1} \ \frac{q^{-\frac{1}{2}(i-n-2)(i-n-1)}(1+a'q)}{q^{-\frac{1}{2}
                       [(i-1)i+n(n+1)]}(1+aq)} \, f^{(i-n-1)}m_0 \
                      \mbox{for some} \ a' \in \ba\\
                     &=a^*q^{(i-1)(n+1)} f^{(i-n-1)}m_0 \ \mbox{for some} \ a^* \in \ba^\times.
\end{align*}
\textit{Case 2:} \ $0 \le i \le n$
\begin{align}\label{cc16}
 f^{(i)}f^{(-n-1)}m_0&=(-1)^{n+1}(-1)^{n+1-i} \ \frac{q^{-\frac{1}{2}(n-i)(n+1-i)}(1+a'q)}
                       {q^{-\frac{1}{2}[(i-1)i+n(n+1)]}(1+aq)} \, f^{(i-n-1)}m_0 \
                      \mbox{for some} \ a' \in \ba\notag\\
                     &=a^*q^{ni} f^{(i-n-1)}m_0 \ \mbox{for some} \ a^* \in \ba^\times.
\end{align}
Hence, for $r(q)=q^{-n(n+1)}$, we have by Lemmas \ref{lattices} and \ref{l:*}
\begin{align*}
 \widetilde{L}&=
    \oplus_{i \ge 0} \ba q^{-n(n+1)}f^{(i)}f^{(-n-1)}m_0\\
 &=\left(\oplus_{i=0}^n \ba q^{n(i-n-1)}f^{(i-n-1)}m_0\right)\oplus\left(\oplus_{i \ge n+1}
     \ba q^{(n+1)(i-n-1)}f^{(i-n-1)}m_0\right)\\
 &=\oplus_{i \ge 0} \ba \frac{q^{n(i-n-1)}}{1+q^{-(i-n-1)}} \, f^{(i-n-1)}m_0
  =\oplus_{i \ge 0} \ba \frac{2q^{n(i-n-1)}}{1+q^{-(i-n-1)}} \, f^{(i-n-1)}m_0.
 \end{align*}

\item Evidently, $L^e=\ba m_0 \subseteq \w{L}$, $\w{L} \ns L$, and
$\w{L} \cap M=\oplus_{j \ge 0} \ba q^{j(n+1)}f^{(j)}m_0=\w{L} \cap L$.
Next,
\begin{align*}
 \w{L} \slash (\w{L} \cap L) &=\oplus_{i=0}^n \ba q^{-n(n+1-i)}f^{(-n-1+i)}\overline{m_0}
                               \quad \mbox{where} \ \overline{m_0}=m_0+\w{L}\cap L\\
                             &=\oplus_{i=0}^n \ba q^{-n(n+1)}f^{(i)}f^{(-n-1)}\overline{m_0} \quad \mbox{by} \   (\ref{cc16}).
\end{align*}

Since $q^{-n(n+1)}f^{(-n-1)}\overline{m_0}$ (in fact, its image under $\w{L} \slash
(\w{L} \cap L) \hookrightarrow C(M) \slash M$)
is a highest weight vector of $C(M) \slash M$, then $\w{L} \slash (\w{L} \cap L)$
is a crystal lattice of $C(M) \slash M$. Thus, $\w{L}$ is a completion of $L$.

\end{enumerit}
\end{proof}

\begin{rem}
\begin{enumerit}
\item The above proposition gives the basis vectors of a crystal lattice of $C(M)$ in a closed form.
However, if we remove the $\ba^\times$-multiples from the given basis vectors of $\w L$, they come down to
\begin{equation}\label{cc17}
 \begin{cases}
  q^{n(i-n-1)}f^{(i-n-1)}m_0       &\mbox{for} \ 0 \le i \le n\\
  q^{(n+1)(i-n-1)}f^{(i-n-1)}m_0   &\mbox{for} \ i > n
 \end{cases}
\end{equation}
and they generate the same crystal lattice $\widetilde L$.

\item If we start with the lattice from Proposition \ref{p:neg} and apply $f^{(n+1)}$
to its highest weight vector, we get:

 $$f^{(n+1)}\left(\frac{2q^{-n(n+1)}}{1+q^{n+1}} \, f^{(-n-1)}m_0\right)
                  =\frac{2q^{-n(n+1)}}{1+q^{n+1}} \frac{[0]!}{[n+1]![-n-1]!} \, m_0
                 =a^*m_0 \ \mbox{for some} \ a^* \in \ba^\times.$$
In other words, $m_0 \in \w{L}.$

\end{enumerit}
\end{rem}

%%%%%%%%%%%%%%%%%%%%%%%%%%%%%%%%%%%%%%%%%%%%%%%

\subl{}

In the last subsection, we show how an operator defined by Kashiwara in \cite{K1}
can be applied to a certain $\ba$-lattice in order to deform it to a crystal
lattice of the completion.

We consider the $\ba$-lattices $L$ in $M=M(-n-2)$ and $L^\sharp$ in
$C(M)$ where $L=\oplus_{k \ge 0} \ba f^{(k)}m_0$ and
$L^\sharp=\oplus_{k \ge -n-1} \ba f^{(k)}m_0$. Note that $L$ is a
crystal lattice of $M$, but $L^\sharp$ is \emph{not} its completion.
Actually, $L^\sharp$ is not even a crystal lattice since
$\w{e}(L^\sharp) \nsubseteq L^\sharp$ where $\w{e}$ denotes the
Kashiwara operator on $C(M)$. However, we will directly transform
$L^\sharp$ to a completion of $L$.

Let $\triangle = qt+q^{-1}t^{-1}+(q-q^{-1})^2fe-2$ be a central element of $U_q$
and consider the action of $qt\triangle$ on $C(M)$.
For $k>0$,
$f^{k}\triangle=\triangle f^{k}$, thus $\triangle
f^{-k}m_0=f^{-k}\triangle m_0$.
Hence $\triangle f^{(k)}m_0=f^{(k)}\triangle m_0
=\left(q^{-n-1}+q^{n+1}-2\right)f^{(k)}m_0$ for $k\ge -n-1$.

Now, $qt\triangle f^{(k)}m_0=\left(q^{-n-1}+q^{n+1}-2\right)q^{-n-1-2k}f^{(k)}m_0
                       =\left[(q^{-n-1}-1)q^{-k}\right]^2 f^{(k)}m_0$.
We define $(qt\triangle)^{\frac{1}{2}} \in \mbox{End} \, \bigl(C(M)\bigr)$ as
follows: $(qt\triangle)^{\frac{1}{2}}f^{(k)}m_0= q^{-k}(q^{-n-1}-1)
f^{(k)}m_0$ for $k \ge -n-1$ (cf.\cite{K1}).
Thus,
$(qt\triangle)^{-\frac{1}{2}} \in \mbox{End} \, \bigl(C(M)\bigr)$ is defined by
\begin{equation}\label{n1}
 (qt\triangle)^{-\frac{1}{2}}f^{(k)}m_0=q^k(q^{-n-1}-1)^{-1}f^{(k)}m_0
\end{equation}
for $k \ge -n-1$. Let $S_n \in \mbox{End} \, \bigl(C(M)\bigr)$ be given by
\begin{equation}\label{n2}
 S_n=
  \begin{cases}
   q^{-n-1} (qt\triangle)^{-\frac{1}{2}},    &\mbox{on} \ C(M)_\ell \ \, \mbox{for} \, \ell \le -n-2\\
   \mbox{id},                         &\mbox{on} \ C(M)_\ell \ \, \mbox{for} \, \ell > -n-2.
  \end{cases}
\end{equation}

\begin{prop}
Let $L=\oplus_{k \ge 0} \ba f^{(k)}m_0$ be a crystal lattice of the
Verma module $M$ with highest weight $-n-2$, and let $L^\sharp=\oplus_{k
\ge -n-1} \ba f^{(k)}m_0$. Then $\displaystyle{q^{-n(n+1)}S_n
(qt\triangle)^{-\frac{n}{2}} L^\sharp}$ is a completion  of $L$.
\end{prop}
\begin{proof}
We note that $q^{-n-1}-1=a q^{-n-1}$ where $a=1-q^{n+1} \in
\ba^\times$. It follows from (\ref{n1}) that
$(qt\triangle)^{-\frac{n}{2}} L^\sharp =\oplus_{k \ge -n-1} \ba
q^{nk}(q^{-n-1})^{-n} f^{(k)}m_0=\oplus_{k \ge -n-1} \ba
q^{n(k+n+1)}f^{(k)}m_0$. Also, for $k \ge 0$,
\begin{equation}\label{n5}
 S_n f^{(k)}m_0=q^{-n-1}q^k
                    (q^{-n-1}-1)^{-1}f^{(k)}m_0=q^{k}a^{-1}f^{(k)}m_0.
\end{equation}
Therefore,
\begin{align*}
 q^{-n(n+1)}S_n (qt\triangle)^{-\frac{n}{2}}
 L^\sharp &=S_n\left(\oplus_{k \ge -n-1} \ba q^{nk}f^{(k)}m_0\right)\\
          &=\left(\oplus_{k \ge 0} \ba q^{(n+1)k}f^{(k)}m_0\right) \oplus
            \left(\oplus_{-n-1 \le k < 0} \ba q^{nk}f^{(k)}m_0\right).
\end{align*}
Hence, $q^{-n(n+1)}S_n (qt\triangle)^{-\frac{n}{2}} L^\sharp$ is the $\ba$-lattice generated by
the same vectors as in (\ref{cc17}). Therefore, it is equal to the crystal lattice $\w{L}$
from the Proposition \ref{p:neg} which is a completion of $L$.
\end{proof}

\vspace{.3in}

{\bf Acknowledgments:} \ My sincere thanks are due to V. Deodhar for
many inspiring and useful discussions, and in particular for his
contribution to Section~2, and to A. Moura for reading through the
manuscript and for valuable suggestions. The referee's comments led
to a clearer exposition and are appreciated. I am also pleased to
thank the Max-Planck-Institut f\"ur Mathematik in Bonn for its
hospitality.

\vspace {.5in}

%%%%%%%%%%%%%%%%%%%%%%%%%%%%%%%%%%references%%%%%%%%%%%%%%%%%%%%%%%%%%%%%

\end{document}